\documentclass[10pt]{article}
\usepackage[b5paper]{geometry}
\usepackage{amsmath,amsthm,amsfonts,amssymb,MnSymbol,mathrsfs,latexsym}
\usepackage[english]{babel}

\usepackage[applemac]{inputenc}
\usepackage[all]{xy}
\numberwithin{equation}{section}

\DeclareMathOperator*{\esssup}{ess\,sup}

  \newcommand{\const}{\rm const}

\theoremstyle{plain}
\newtheorem{theorem}{Theorem}[section]
\theoremstyle{theorem}

\newtheorem{remark}{Remark}[section]

\renewenvironment{proof}{{\bf{Proof.}}}{\hfill $\Box$ \\}

\begin{document}

\title{ \large \textbf{Grand Lebesgue Spaces are really Banach algebras relative to the convolution on unimodular locally compact groups}}

\footnotesize\date{}

\author{\normalsize \textbf{Maria Rosaria FORMICA ${}^{1}$},   \normalsize\textbf{Eugeny OSTROVSKY
${}^2$} and \normalsize\textbf{Leonid SIROTA ${}^3$}}

\maketitle

\begin{center}
{\footnotesize ${}^{1}$ Universit\`{a} degli Studi di Napoli \lq\lq Parthenope\rq\rq, via Generale Parisi 13,\\
Palazzo Pacanowsky, 80132,
Napoli, Italy.} \\

\vspace{2mm}

{\footnotesize e-mail: mara.formica@uniparthenope.it} \\

\vspace{4mm}

{\footnotesize ${}^{2,\, 3}$  Bar-Ilan University, Department of Mathematics and Statistics, \\
52900, Ramat Gan, Israel.} \\

\vspace{2mm}

{\footnotesize e-mail: eugostrovsky@list.ru}\\

\vspace{2mm}

{\footnotesize e-mail: sirota3@bezeqint.net} \\

\vspace{4mm}

\end{center}

\begin{abstract}
We prove that the Grand Lebesgue Space, builded on a unimodular
locally compact topological group, forms a Banach algebra relative
to the convolution.
\end{abstract}


\vspace{4mm}


\vspace{2mm} \noindent {\it 2010 Mathematics Subject
Classification}:
46E30,   
44A35,     
43A10.    

\vspace{5mm}

\section{Definitions. Statement of problem. Notations.}

\vspace{5mm}

Let $G$ be a unimodular locally compact topological group equipped
with bi-invariant Borelian Haar measure $\mu$. Define, as usually,
the convolution between two measurable integrable functions $ \
f,g:G \to R \ $ by
$$
(f*g)(x) \stackrel{def}{=} \int_G f(y) g(y^{-1} \, x) \, d\mu(y).
$$

The usual Lebesgue-Riesz spaces $L_p = L_p(G)=L_p(G,\mu)$ are the
spaces of all measurable functions $f: G\to \mathbb R$ defined by
the norms
$$
 \begin{array}{ll}
      |f|_p := \displaystyle \left[ \int_G |f(x)|^p
\, d\mu(x) \right]^{1/p}, & 1\leq p<\infty \\ \\
|f|_{\infty} := \displaystyle \esssup_{x \in G} |f(x)|, & p=\infty
    \end{array}
$$

The well-known Young (or Hausdorff-Young) inequality has the form

\begin{equation} \label{HausYoung}
|f*g|_r \le |f|_p \ |g|_q, \ 1 \le p,q,r \le \infty,
\end{equation}
where
\begin{equation} \label{restiction}
 1+1/r = 1/p + 1/q,
\end{equation}
see e.g. \cite{Hardy}, \cite{Zelazko}, \cite{Bennett Sharpley1988}
and \cite{Beckner}, \cite{Brascamp}, \cite{Fournier} for the
computation of the best possible constant in (\ref{HausYoung}).\par

\noindent Generalized convolution inequalities and applications have
been obtained in \cite{Jain-Jain-2017} and fractional convolutions
in \cite{Jain-Jain-Kumar-2015}.

Choosing $r = p, \, q = 1$ in (\ref{HausYoung}), we deduce
\begin{equation} \label{HausYoungcomp}
|f*g|_p \le |f|_p \ |g|_1.
\end{equation}

If in addition the group $G$ is also compact, one can suppose
$\mu(G)
= 1$ so that 
$|g|_1 \le |g|_p$. Therefore
\begin{equation} \label{HausYoungBanAlg}
|f*g|_p \le |f|_p \ |g|_p, \  p \in [1,\infty].
\end{equation}

The last inequality shows that the space $L_p(G,\mu)$, relative to
the convolution bilinear operator, forms in this case a {\it Banach
algebra}, see \cite{Zelazko}.

\vspace{2mm}

Our aim is to extend inequality \eqref{HausYoungBanAlg} to the
so-called Grand Lebesgue Spaces, improving some results recently
obtained in \cite{Turan}; see also \cite{Gurkanli1},
\cite{Gurkanli2}. \par

 \vspace{2mm}

In more details, we will find a sufficient condition on the
generating function of these spaces under which they form also a
 {\it Banach algebra}, relative to the convolution. \par

\begin{center}

 \vspace{2mm}
\ {\sc Brief note about  Grand Lebesgue Spaces. } \\

\end{center}

\vspace{2mm}

 We recall here, for reader convenience, some known definitions and facts about the theory of the Grand Lebesgue Spaces which we will use.
 We consider only the case of the space on the set $G$ equipped  with the {\it Haar} measure $\mu$. \par

 Let $\psi = \psi(p) = \psi[b](p), \ p \in [1,b), \ b = \const \in (1,\infty]$, be a numerical valued non-negative function, not necessarily finite in every point, such that

\begin{equation} \label{Positive psi}
 \inf_{p \in [1,b)} \psi(p) > 0
\end{equation}
and

\begin{equation} \label{Positive norming}
  \psi(1) \in( 0, \ \infty).
\end{equation}
The set of all such functions will be denoted by $ \Psi = \Psi[b] =
\{\psi(\cdot) \}$. If $b < \infty$, one can assume $p \in [1,b]$.
\par

\vspace{2mm}

We will assume, without loss of generality, \, $\psi(1) = 1$ for all
such functions $\psi(\cdot)$.

 \vspace{2mm}

\noindent For instance, the following functions
$$
\psi_{(m)}(p) := p^{1/m}, \ \ m > 0, \  \ p \in [1,\infty)
$$
or
$$
\psi^{(b; \beta)}(p) := (b-1)^{\beta} (b-p)^{-\beta}, \ \ p \in
[1,b), \  \ \beta \ge 0, \ \ b\in (1,\infty)
$$
satisfy \eqref{Positive psi} and \eqref{Positive norming}, with
$\psi(1) = 1$.

\vspace{2mm}

The (Banach) Grand Lebesgue Space $G\psi =G\psi [b]$ consists of all
the real (or complex) numerical valued measurable functions $f: G
\to R$ having finite norm, defined by

 \begin{equation} \label{norm psi}
    || f ||_{G\psi} = ||f||_{G\psi[b]} \stackrel{def}{=} \sup_{p \in [1,b)} \left[ \frac{|f|_p}{\psi(p)} \right].
 \end{equation}

 \vspace{4mm}

 The function $\psi = \psi(p)=\psi[b](p) $ is named as the {\it generating function} for this space. \par

If for instance
$$
 \psi(p) = \psi^{(r)}(p)=\left\{
  \begin{array}{ll}
    1  & , \ \ p=r \\
    +\infty &  , \ \ p\neq r
  \end{array}
\right.
$$
where $ r\in [1,\infty),  \ C/\infty := 0, \ C \in \mathbb R$,
(extremal case), then the corresponding $ G\psi^{(r)}(p)$ space
coincides with the classical Lebesgue-Riesz space $L_r = L_r(G)$.

\vspace{4mm}

Note that the space $L^{b),\theta}$ introduced in
\cite{Greco-Iwaniec-Sbordone-1997} (see \cite{Iwaniec-Sbordone-1992}
for $\theta =1$), through the norm

$$
||f||_{L^{b),\theta}}=||f||_{b, \theta,G} \stackrel{def}{=} \sup_{0
< \epsilon \le b - 1} \left[ \ \epsilon^{\theta/(b-\epsilon)} |f|_{b
- \epsilon} \ \right], \ \theta \ge 0
$$
quite coincides, up to equivalence of the norms, with the
appropriate one $G\psi^{(b; \beta)}$, with $\beta = \theta/b$.

Namely, 
$$
C^-(b,\theta,G)||f||_{G\psi^{(b; \beta)}}\leq ||f||_{L^{b),\theta}}
\leq C^+(b,\theta,G)||f||_{G\psi^{(b; \beta)}},
$$
where $0<C^-(b,\theta,G)\leq C^+(b,\theta,G)<\infty$.
 \vspace{4mm}

Let $f: G \to R$ be a measurable function and $b \in (1,\infty]$
such that
$$
|f|_p < \infty,  \ \ \forall p \in [1,b);
$$
the so-called {\it natural function} $\psi^{(f)}(p)$ for $f$ is
defined by
$$
\psi^{(f)}(p) \stackrel{def}{=} |f|_p.
$$
Obviously,
$$
||f||_{G\psi^{(f)}} = 1.
$$

\vspace{4mm}

 The Grand Lebesgue Spaces have been widely investigated, e.g.
\cite{Fiorenza-2000-duality}, \cite{Capone-Fiorenza-2005-small},
\cite{Capone-Formica-Giova-2013}, \cite{Formica-Giova-2015-Boyd},
\cite{Iwaniec-Koskela_Onninen-2001}, \cite{KozOs},
\cite{LiflOstSir}, \cite{Ostrovsky1} - \cite{Ostrovsky5},
\cite{Jain-Singh-Singh_Stepanov-2019}, \cite{Molchanova2019}, etc.
They play an important role in the theory of Partial Differential
Equations (PDEs) (see e.g. \cite{Fiorenza-Mercaldo-Rakotoson-2001},
\cite{Fiorenza-Mercaldo-Rakotoson-2002},
\cite{Greco-Iwaniec-Sbordone-1997}, \cite{fioformicarakodie2017}),
in interpolation theory (see e.g. \cite{Fiorenza-Karadzhov-2004},
\cite{Fiorenza-Formica-Nonlinear_Anal_2018}), in the theory of
Probability (\cite{Ermakov}, \cite{KozOsSir2017},
\cite{Ostrovsky3}-\cite{Ostrovsky5}), in Statistics and in the
theory of random fields (see e.g. \cite{KozOs}, \cite[chapter
5]{Ostrovsky1}), in Functional Analysis and so on.
\par

 These spaces are rearrangement invariant (r.i.) Banach functional spaces; its fundamental function has been studied in  \cite{Ostrovsky4}. They
 not coincides, in the general case, with the classical rearrangement invariant spaces: Orlicz, Lorentz, Marcinkiewicz, etc., see \cite{LiflOstSir}, \cite{Ostrovsky2}.

\vspace{5mm}

\section{Main result.}

\begin{theorem}\label{main result}

Let $G$ be a unimodular locally compact topological group equipped
with bi-invariant Haar measure $\mu$ and let $G\psi $ be the
(Banach) Grand Lebesgue space builded on $G$. If $ \psi(1) = 1$,
then $G\psi$ is also a Banach algebra under the convolution as a
multiplicative operation, that is
\begin{equation} \label{result}
||f * g||_{G\psi} \le ||f||_{G\psi} \cdot ||g||_{G\psi}
\end{equation}
for all $f,g: G\to \mathbb R$ in $G_{\psi}$.
\end{theorem}

\vspace{4mm}

\begin{proof}
 Let $f,g$ be two non-zero functions from the space $ G\psi$;  we assume, without loss of generality, that
$$
||f||_{G\psi} = ||g||_{G\psi} = 1.
$$
From the definition of the norm in the Grand Lebesgue Spaces it
follows immediately that
$$
|f|_p \le \psi(p), \ \ |g|_p \le \psi(p), \  \ p \in [1,b).
$$
We apply  once again the Young inequality (\ref{HausYoungcomp})
$$
|f*g|_p \le |f|_p \, |g|_1 \le \psi(p) \, |g|_1,
$$
thus
$$
\frac{|f*g|_p}{\psi(p)} \le |g|_1 \le \psi(1) = 1
$$
and, taking the supremum over $p \in [1,b)$, we have

$$
\sup_{p \in [1,b)} \left\{ \ \frac{|f*g|_p}{\psi(p)} \ \right\}
\le |g|_1 \le \psi(1) = 1.
$$

Therefore

\begin{equation*} 
||f*g||_{G\psi} \le 1 = ||f||_{G\psi} \cdot ||g||_{G\psi}.
\end{equation*}

\end{proof}

\vspace{2mm}

\begin{remark}
 In the next section we will see that equality in \eqref{result} can be achieved.
\end{remark}

\section{An example.}

\vspace{5mm}

 Let now  $ G =\mathbb R$ with the classical Lebesgue measure $d \mu = dx$ and the ordinary convolution
\begin{equation}\label{convolution in R}
(f*g)(x)=\int_{-\infty}^{+\infty} f(y)g(x-y)\, dy,
\end{equation}
for all measurable integrable functions $f,g$ defined on $\mathbb
R$.

Define the Gaussian  density
\begin{equation}\label{gaussian_density}
z =z(x) = z_{\sigma}(x) := (\sigma \ \sqrt{2 \ \pi})^{-1} \, \exp
\left\{ \frac{-x^2}{2\sigma^2}  \right\}, \ \ \sigma = \const > 0;
\end{equation}
then its natural function $\psi^{(z,\sigma)}(p)$ has the form

\begin{equation}\label{gaussian_psi}
\psi_{\sigma}(p)= \psi^{(z_{\sigma})}(p) = | \ z_{\sigma}(\cdot) \
|_p = (2 \pi)^{ 1/(2p) - 1/2 } \ p^{-1/(2p)} \sigma^{1/p - 1}, \ \ p
\in [1,\infty).
\end{equation}

One can define formally the value $\psi^{(z_{\sigma})}(p)$ for $p
= \infty$ as

$$
\psi_{\sigma}(\infty)  \stackrel{def}{=} \lim_{p \to \infty} \psi_{\sigma}(p) = (\sigma \ \sqrt{2 \ \pi})^{-1} = \max_{x \in R} z_{\sigma}(x).
$$

 Notice that $\psi_{\sigma}(1) = 1$. By virtue of Theorem \ref{main result} we conclude that the Grand Lebesgue Space $G\psi^{(z_{\sigma})}(p)$, builded over
 the ordinary real line $ \mathbb R$, relative to the convolution operation, forms also a Banach algebra, despite the (commutative) group $(\mathbb R,+)$ is not compact.

\vspace{2mm}

Moreover, let us choose $f(x) = g(x) = z_{1}(x)$, then

$$
||f||_{G\psi_{1}} = ||g||_{G\psi_{1}} = 1.
$$
Furthermore,
$$
h(x) := (f*g)(x)= z_{\sqrt{2}}(x);
$$

$$
|h|_p = (2 \pi)^{ 1/(2p) - 1/2 } \ p^{-1/(2p)} [\sqrt{2}]^{1/p - 1}, \ p \in [1,\infty),
$$
so that

$$
  \frac{|h|_p}{\psi_{1}(p)} = [\sqrt{2}]^{1/p - 1},
$$
and

$$
||h||_{G\psi_{1}}  = \sup_{p \ge 1} [\sqrt{2}]^{1/p - 1} = 1.
$$
 Thus, in this case, we have

$$
||h||_{G\psi_{1}} = ||f*g||_{G\psi_{1}}  = ||f||_{G\psi_{1}} \cdot
||g||_{G\psi_{1}}.
$$

\vspace{2mm}

 Note that the case $ \ G = \mathbb R^d \ $ may be investigated analogously.

\vspace{5mm}

\section{Concluding remarks.}

\vspace{5mm}

 \ {\bf A.} Let $E$ be an arbitrary translation invariant Banach functional space  with norm $||\cdot||_E $ builded over $ \ G \ $
 and such that

$$
\forall g \in L_1 \ \Rightarrow |g|_1 \le ||g||_E.
$$
 We have, using the triangle inequality,

\begin{equation*}
\begin{split}
||f*g||_E = & || \int_G f(y^{-1} x) \ g(y) \ \mu(dy)  ||_E \le
\int_G
||f||_E \ |g(y)| \ \mu(dy) \\
= &||f||_E \ |g|_1 \le    ||f||_E \cdot ||g||_E.
 \end{split}
\end{equation*}
 Therefore, $E$ is again a Banach algebra.

\vspace{4mm}

{\bf B.} Let us prove the {\it necessity} of the condition
\eqref{restiction} at least in the case $G= \mathbb R^d, \ d =
1,2,3,\ldots$. Suppose that for {\it arbitrary} non-zero and
non-negative numerical valued smooth functions $f,g$, defined on
$\mathbb R^d$ and having compact support, it is

\begin{equation} \label{estimat}
|f*g|_r \le |f|_p \ |g|_q, \ \   p,q,r \ge 1.
\end{equation}
Then $ \ 1+ 1/r = 1/p + 1/q$. \par

Actually, one can use the so-called scaling method, see e.g.
\cite{Talenti}. Assume that (\ref{estimat}) is satisfied. Define for
any positive constant $ \lambda > 0 $ the dilation operator

$$
T_{\lambda}[f](x) = f(\lambda \ x),  \ \ x \in \mathbb R^d.
$$
We deduce from \eqref{estimat}

\begin{equation} \label{estim}
|T_{\lambda}f*T_{\lambda}g|_r \le |T_{\lambda}f|_p \
|T_{\lambda}g|_q, \  \  p,q,r \ge 1
\end{equation}
and, after simple calculations,

$$
\lambda^{-d - d/r} \ |f*g|_r \le \lambda^{-d/p - d/q} \ |f|_p \ |g|_q.
$$
Since $\lambda$ is a arbitrary positive value, we conclude $ \  -d -
d/r = -d/p - d/q,  \ $ which is quite equally to
\eqref{restiction}.\par

\vspace{5mm}

 {\bf C.} As long as the Grand Lebesgue Space  $ \ G\psi_{1}$ forms a Banach algebra, it is interesting,
by our opinion, to consider its ideals, maximal or not. We apply
here well-known classical results from \cite{Cigler}, \cite{Hewitt},
\cite{Hewitt1970}, \cite{Naimark}.\par

 Define as usually for this purpose the Fourier transform for every function from this space

\begin{equation} \label{Four}
\hat{f}(\xi) = F[f](\xi) := \int_{\mathbb R^d} e^{ 2 \pi i \ (\xi,x)
} \ f(x) \  dx,
\end{equation}
where $  \  (\xi,x) \ $ denotes the inner (scalar) product and $ \  |x| := \sqrt{(x,x)}.  \ $ \par

Since  $f \in L_1(\mathbb R^d)$, the integral in \eqref{Four}
converges uniformly and it represents a uniformly continuous bounded
function on the whole space $\mathbb R^d$

$$
\sup_{\xi \in \mathbb R^d}|F[f](\xi)| \le |f|_1 \le ||f||G\psi_1 <
\infty.
$$

Here $\psi_1 = \psi_1(p), \ p \in [1,\infty)$ is the function
defined by \eqref{gaussian_psi}, with $\sigma=1$. For this purpose
one can replace $ \ \psi_1(\cdot) \ $ with many other such
functions, for instance, with the more general $ \
\psi_{\sigma}(\cdot), \ \sigma \in (0,\infty). \ $
\par

Recall that if $h(x) = (f * g)(x) $, then

\begin{equation} \label{FourConv}
F[h](\xi) = F[f](\xi) \cdot F[g](\xi).
\end{equation}

Let $\eta$ be an arbitrary fixed point of the  space $\mathbb R^d$;
define the set

\begin{equation} \label{IdealOrd}
J(\eta) = \{ f \,:\, f \in G\psi_1, \ F[f](\eta) = 0  \}.
\end{equation}

It follows immediately from the  relation \eqref{FourConv} that the
set $J(\eta)$ forms a well-known classical (maximal) ideal in the
Banach algebra $G\psi_1$. \par

\vspace{4mm}

Denote by $ \ UC = UC(\mathbb R^d) \ $ the set of all {\it
uniformly} continuous bounded functions $ \ f: \mathbb R^d \to
\mathbb R \ $ or $ \ f: \mathbb R^d \to \mathbb C. \ $ This set
forms also an ideal in this algebra. Namely, define the usual
modulus of continuity for a function $ \ f \in UC(\mathbb R^d) \ $

$$
\omega[f](\delta) \stackrel{def}{=} \sup_{|x-y|\le \delta} |f(x) -
f(y)|, \  \ \delta \in [0,\infty)
$$
so that
$$
f \in  UC(\mathbb R^d) \ \Rightarrow \ \lim_{\delta \to 0+}
\omega[f](\delta) = 0.
$$
For $ f \in  UC(\mathbb R^d)$ the function $x \to h(x) = (g*f)(x)$
is bounded and uniformly continuous, in fact if
$$
x_1,x_2 \in \mathbb R^d \ : \   |x_1 - x_2|\le \delta
$$
we have
\begin{equation*}
\begin{split}
|h(x_1) - h(x_2)| \le & \int_{\mathbb R^d} |g(y)| |f(x_1 - y) -
f(x_2 - y)| \, dy \\
\le & |g|_1 \omega[f](\delta) \le ||g||_{G\psi_1}  \cdot
\omega[f](\delta) \to 0, \ \ \delta \to 0+.
\end{split}
\end{equation*}
Thus, $h(\cdot) \in UC(\mathbb R^d)$.

\vspace{5mm}

{\bf D.} Let us show  now that the condition \eqref{Positive
norming} is, in the general case, essential. Namely, suppose
temporarily that $ \ \psi(1) = 0, \ $ then

$$
 \forall f \in G\psi \ \Rightarrow \ |f|_1 = 0
$$
and, consequently, $ f = 0$, i.e.  $G\psi = \{ \ 0 \ \}$.
\par
Furthermore, let us consider the following counterexample. Choose $
G = \mathbb R$ with the usual convolution. Define the function

$$
f(x)=f_{3/2}(x):= x^{-2/3} \, I(x \ge 1),
$$
where $I(A)$ denotes the indicator function of the predicate $A$.
\par

We deduce, after simple calculations,
\begin{equation*}
\tilde{\psi}(p) := |f|_p =\left\{
  \begin{array}{ll}
    \left( \displaystyle \frac{3}{2p-3}
\right)^{1/p}, & p > 3/2 \\ \\
    +\infty, & p \in [1,3/2].
  \end{array}
\right.
\end{equation*}
Evidently,
$$
||f||_{G\tilde{\psi}} = 1,
$$
as long as the function $\tilde{\psi}(p)$ is the natural function
for $f$. \par

Let us now evaluate the convolution
$$h(x) = h(x)=(f*f)(x)$$
as $ \ x \to \infty$; of course, $h(x) = 0, \ x < 1$.
 For $x \geq 3$ we have, as $ x \to \infty$,

\begin{equation*}
\begin{split}
h(x) = & \int_1^{x-1}  y^{-2/3} \ (x-y)^{-2/3} \ dy = x^{-1/3}
\int_{1/x} ^{1 - 1/x} \ z^{-2/3} \ (1 - z)^{-2/3} \  dz  \\
\sim & B(1/3, 1/3) \ x^{-1/3} = C_1 \ x^{-1/3},
\end{split}
\end{equation*}
where $B(\cdot, \cdot)$ denotes the ordinary Beta-function and $ \
C_1 = B(1/3, 1/3) \in (0,\infty)$.  \par

More precisely note that, as $x \ge 3$,

\begin{equation*}
h(x) \ge x^{-1/3} \int_{1/3}^{2/3} z^{-2/3} \ (1 - z)^{-2/3} \ dz  =
C_2 \ x^{-1/3}, \ \  C_2 = \const \in (0,\infty);
\end{equation*}

$$
|h|_p^p \ge C_2^p \int_{1/3}^{2/3} x^{-p/3} \ dx =  C_2^p
\frac{3}{p-3} \cdot \frac{1 - 2^{1 - p/3}}{3^{1 - p/3}} =:
 K^p(p), \ \ \ K(p)\in (0,+\infty);
$$

$$
|h|_p \ge K(p), \  \ p > 3,
$$

$$
 |h|_p = \infty,  \ 1 \le p \le 3.
$$

 It remains to note that the function $ \ h(\cdot) \ $ does not belongs to the Grand Lebesgue Space  $ G\tilde{\psi}$. In fact,

$$
||h||_{G\tilde{\psi}} = \sup_{p > 3/2} \left\{
\frac{|h|_p}{\tilde{\psi}(p)} \ \right\}   \ge
\frac{|h|_2}{\tilde{\psi}(2)} = \infty,
$$
or equivalently
 $$
\infty = ||f*f||_{G\tilde{\psi}} > ||f||_{G\tilde{\psi}} \cdot
||f||_{G\tilde{\psi}} = 1.
$$

Thus, the (Banach) Grand Lebesgue Space $ \ G\tilde{\psi} \ $ is not
a Banach algebra. \par

The case when
$$
f(x) = f_{\alpha}(x) = x^{-1/\alpha} \ I(x \ge 1), \ \ \alpha =
\const \in (1,2)
$$
may be investigated quite analogously. \par

By the way note that, despite $f_{\alpha}(\cdot) \notin L_1(R)$, its
Fourier transform $ F[f_{\alpha}](\xi)$ there exists as an element
of the space $L_2(R)$ because $f_{\alpha}(\cdot) \in L_2(R)$, by the
Theorem of  Plancherel-Parseval. Furthermore,
 $ F[f_{\alpha}](\xi) $ there exists for all non-zero values of the variable $ \xi$. More exactly, as $ \xi \to 0+ $

$$
F[f_{\alpha}](\xi) \sim \xi^{1/\alpha - 1}\cdot \int_0^\infty
e^{2\pi i y} y^{-1/\alpha}\,dy
$$
and

$$
\xi \to \infty \ \Rightarrow \ F[f_{\alpha}](\xi) \sim \frac{e^{2
\pi i \xi}}{2 \pi i \xi}.
$$

\vspace{6mm}

 {\it An open and interesting problem, by our opinion, is to find the necessary conditions imposed on the generating function
 $ \ \psi(\cdot) \ $  under which $ \ G\psi  = G\psi_G \ $  forms a Banach algebra and to find, in this
case, all its  maximal ideals.} \par

\vspace{6mm}

\emph{Acknowledgement.} {\footnotesize The first author has been
partially supported by the Gruppo Nazionale per l'Analisi
Matematica, la Probabilit\`a e le loro Applicazioni (GNAMPA) of the
Istituto Nazionale di Alta Matematica (INdAM) and by Universit\`a
degli Studi di Napoli Parthenope through the project \lq\lq sostegno
alla Ricerca individuale\rq\rq .\par

\end{document}